\documentclass{ip-journal}

\begin{document}

\title{Three-Diffeomorphism Conformal Space over Lorentzian Manifold}
\author{Lukasz Andrzej Glinka}

\address{}
\email{laglinka@gmail.com}

\date{\today}

\begin{abstract}
Through making use of a Borel measure and a piecewise-Riemannian inner scalar product, it is shown that over a Lorentzian manifold every three diffeomorphisms generate a conformal space, whose elements are smooth vector-valued functions equipped with compact supports. Few examples, in particular a diffeoinvariant measure, are provided with respect to an arbitrary smooth function introduced as into consideration as a multiplier to a local scale factor.
\end{abstract}

\maketitle
\section{Introduction}

Let us consider a Lorentzian manifold $\mathcal{M}$ of $\mathrm{D}=\mathrm{d}+1$ dimensions and the Sylvester signature $(\mathrm{d},1,0)$, or $\mathrm{s}=\mathrm{d}-1$. Let $g_{\mu\nu}(x)$ be a metric tensor on $\mathcal{M}$ of dimension $\mathrm{D}\times \mathrm{D}$ and properties $g(x)=\det g_{\mu\nu}(x)<0$, $g^\mu_\alpha=\delta^\mu_\alpha$, $g^\mu_\mu=\mathrm{D}$, where $\delta^\mu_\alpha=\mathrm{diag}[1,\ldots,1]$ is the Kronecker $\mathrm{D}\times\mathrm{D}$ unit square matrix. In analysis of such an object one can make use of the theory of function spaces, particularly the class of vector-valued Banach functional spaces, Cf. e.g. the Refs. \cite{b1}-\cite{b3}, known for some modern applications, Cf. e.g. the Refs. \cite{t1}-\cite{t5}.

In this article, $\mathcal{M}$ is concisely approached by the topological techniques of a Borel measure of integration and a piecewise-Riemannian inner scalar product enriched by an arbitrary smooth function which is a multiplier to a local scale function. As a result, any three diffeomorphisms are done the reason to a local scale factor which generates a conformal space $\mathcal{C}(\mathcal{M})$ over $\mathcal{M}$ whose elements are smooth vector-valued functions equipped with compact supports. Furthermore, as the particular case, it is shown that a diffeoinvariant measure over $\mathcal{C}(\mathcal{M})$ can be reached with a help of a certain special choice of an arbitrary smooth function, or a local scale function. Additionally, basic topological, set-theoretical, and analytic properties of $\mathcal{C}(\mathcal{M})$ are briefly discussed.

\section{Borel measure}

Let us consider a coordinate chart $x^\mu\in\mathbb{R}^\mathrm{D}$ on $\mathcal{M}$, and a diffeomorphism $x^\mu\rightarrow \tilde{x}^\alpha(x)$ into a new non-singular coordinate chart $\tilde{x}^\alpha(x)$. Then, a squared arc length on $\mathcal{M}$ is called diffeoinvariant if
\begin{equation}
ds^2=g_{\mu\nu}(x)dx^\mu dx^\nu=\tilde{g}_{\alpha\beta}(\tilde{x})d\tilde{x}^\alpha d\tilde{x}^\beta,\label{jac}
\end{equation}
where usual summation convention holds, and one gets well-known relations
\begin{subequations}
\begin{align}
g_{\mu\nu}(x)&=\tilde{g}_{\alpha\beta}(\tilde{x})\dfrac{\partial\tilde{x}^\alpha}{\partial x^\mu}\dfrac{\partial\tilde{x}^\beta}{\partial x^\nu},\\
g(x)&=\tilde{g}(\tilde{x})\left[\det\dfrac{\partial\tilde{x}^\alpha(x)}{\partial x^\mu}\right]^2.
\end{align}
\end{subequations}

Considering the generic Borel measure in a coordinate chart $\tilde{x}^\alpha(x)$
\begin{equation}
d^{\mathrm{D}}\tilde{x}|\tilde{g}(x)|^{1/2}=dx^0\wedge\ldots\wedge dx^\mathrm{d}|g(x)|^{1/2}\left|\det\dfrac{\partial\tilde{x}^\alpha(x)}{\partial x^\mu}\right|^{-1},\label{meas}
\end{equation}
and performing a change of variables $x^\mu\rightarrow X^\alpha(x)$ such that infinitesimally
\begin{equation}
dx^\mu=\lambda(x)dX^\mu(x),
\end{equation}
where $\lambda(x)\neq0$ is a local scale factor function, one obtains a new Borel measure
\begin{equation}
\mathcal{D}X(x)=dX^0(x)\wedge\ldots\wedge dX^\mathrm{d}(x)[f(x)]^{1/2},\label{meas1}
\end{equation}
where an arbitrary smooth function $f(x)$ is of the form
\begin{equation}
f(x)=[\lambda(x)]^{2\mathrm{D}}|g(x)|\left|\det\dfrac{\partial\tilde{x}^\alpha(x)}{\partial x^\mu}\right|^{-2},
\end{equation}
and conversely
\begin{equation}
\lambda(x)=\left[\dfrac{f(x)}{|g(x)|}\left|\det\dfrac{\partial\tilde{x}^\alpha(x)}{\partial x^\mu}\right|^{2}\right]^{1/(2\mathrm{D})}.
\end{equation}
is determined up to an arbitrary smooth function $f(x)>0$.

This construction shows a non-trivial role of a distinguished collection of $\mathrm{D}$ smooth functions which are components of a smooth vector field
\begin{equation}
X^\mu(x)=\int^x_{x_0}[\lambda(x')]^{-1}dx'^\mu,\label{bomb}
\end{equation}
where $(x_0)^\mu$ is an arbitrary reference point on a $\mathcal{M}$, and $x'$ stands for a variable of integration. Whenever one takes into account the standard definition of support as a closure of the set of non-zero values of a function, then it becomes obvious that all these components are equipped with individual compact supports, and, therefore, the vector-valued smooth function $X^\mu(x)$ in itself is equipped with a compact support.

\section{Conformal space $\mathcal{C}(\mathcal{M})$}

Naturally, a new Borel measure induces a vector-valued function space spanned by a collection of smooth vector fields $\left\{X^\mu(x)\right\}$. Let $Q_{\mu\nu}(X)$ be a metric tensor on this function space. By virtue of diffeoinvariance, a change of variables $x^\mu\rightarrow X^\alpha(x^\mu)$
\begin{equation}
g_{\mu\nu}(x)dx^\mu dx^\nu=Q_{\alpha\beta}(X(x))dX^\alpha(x)dX^\beta(x),
\end{equation}
leads to the relations
\begin{subequations}
\begin{align}
Q_{\alpha\beta}(X(x))&=g_{\mu\nu}(x)\dfrac{\partial x^\mu}{\partial X^\alpha(x)}\dfrac{\partial x^\nu}{\partial X^\beta(x)},\\
Q(X(x))&=g(x)\left[\det\dfrac{\partial x^\mu}{\partial X^\alpha(x)}\right]^2,\label{lob}
\end{align}
\end{subequations}
with $Q(X(x))=\det Q_{\mu\nu}(X(x))$. Meanwhile, for a coordinate chart given by the Eq. (\ref{bomb}) one obtains
\begin{subequations}
  \begin{align}
 \dfrac{\partial x^\alpha}{\partial X^\mu(x)}&=\left[\dfrac{\partial X^\mu(x)}{\partial x^\alpha}\right]^{-1}=\lambda(x)\delta^\alpha_\mu,\label{blo}\\
 \det\dfrac{\partial x^\alpha}{\partial X^\mu(x)}&=[\lambda(x)]^{\mathrm{D}},\label{blo1}
  \end{align}
\end{subequations}
and, consequently
\begin{subequations}
  \begin{align}
  Q_{\alpha\beta}(X(x))&=\left[\lambda(x)\right]^{2}g_{\alpha\beta}(x),\label{kim}\\
  Q(X(x))&=[\lambda(x)]^{2\mathrm{D}}g(x),\label{uop}
  \end{align}
\end{subequations}
what shows that $Q(X(x))<0$ and that $\mathcal{C}(\mathcal{M})$ is a conformal space over $\mathcal{M}$. As a result due to the Eq. (\ref{lob}), one has the Jacobian determinant
\begin{equation}
\left|\det\dfrac{\partial x^\alpha}{\partial X^\mu(x)}\right|=\left|\dfrac{f(x)}{g(x)}\right|^{1/2}\left|\det\dfrac{\partial\tilde{x}^\alpha(x)}{\partial x^\mu}\right|,
\end{equation}
which is compatible with the Eq. (\ref{blo1}). Applying the Jacobian determinant of the composition
\begin{equation}
\det\dfrac{\partial x^\mu}{\partial X^\alpha(x)}\det\dfrac{\partial X^\alpha(x)}{\partial \tilde{x}^\beta(x)}\det\dfrac{\partial \tilde{x}^\beta(x)}{\partial x^\mu}=1,
\end{equation}
one gets the relation
\begin{equation}
\left|\det\dfrac{\partial X^\alpha(x)}{\partial \tilde{x}^\beta(x)}\right|=\left|\dfrac{g(x)}{f(x)}\right|^{1/2}\left|\det\dfrac{\partial \tilde{x}^\beta(x)}{\partial x^\mu}\right|^{-2},
\end{equation}
which leads to a non-diffeoinvariant Borel measure
\begin{equation}
\int\mathcal{D}X|Q(X)|^{1/2}=\int_{\mathcal{M}} d^Dx'|f(x')g(x')|^{1/2}.
\end{equation}

\section{Examples}
For $f(x)=1$, which corresponds to the scale factor
\begin{equation}
\lambda(x)=\left[|g(x)|^{-1}\left|\det\dfrac{\partial\tilde{x}^\alpha(x)}{\partial x^\mu}\right|^{2}\right]^{1/(2\mathrm{D})},
\end{equation}
one obtains the canonical diffeoinvariant Borel measure
\begin{equation}\label{pop}
\int\mathcal{D}X|Q(X)|^{1/2}=\int_{\mathcal{M}} d^Dx'|g(x')|^{1/2},
\end{equation}
along with the formula
\begin{equation}
Q(X(x))=-\left|\det\dfrac{\partial\tilde{x}^\alpha(x)}{\partial x^\mu}\right|^{2}.
\end{equation}
On the other hand, for $f(x)=|g(x)|^{-1/2}$, or
\begin{equation}
\lambda(x)=\left[|g(x)|^{-3/2}\left|\det\dfrac{\partial\tilde{x}^\alpha(x)}{\partial x^\mu}\right|^{2}\right]^{1/(2\mathrm{D})},
\end{equation}
a new Borel measure corresponds with a flat Lorentzian manifold
\begin{equation}\label{pop1}
\int\mathcal{D}X|Q(X)|^{1/2}=\int_{\mathcal{M}} d^Dx',
\end{equation}
whereas
\begin{equation}
Q(X(x))=-|g(x)|^{-1/2}\left|\det\dfrac{\partial\tilde{x}^\alpha(x)}{\partial x^\mu}\right|^{2}.
\end{equation}
For a more sophisticated choice
\begin{equation}
f(x)=|g(x)|\left|\det\dfrac{\partial\tilde{x}^\alpha(x)}{\partial x^\mu}\right|^{-2},
\end{equation}
which gives the scale factor $\lambda(x)=1$, one has identically $Q_{\alpha\beta}(X(x))=g_{\alpha\beta}(x)$ and $Q(X(x))=g(x)$. Because in this case a new Borel measure is
\begin{equation}
\int\mathcal{D}X|Q(X)|^{1/2}=\int_{\mathcal{M}} d^Dx'|g(x')|\left|\det\dfrac{\partial\tilde{x}^\alpha(x')}{\partial x'^\mu}\right|^{-1},
\end{equation}
one sees that only for
\begin{equation}
\left|\det\dfrac{\partial\tilde{x}^\alpha(x)}{\partial x^\mu}\right|=|g(x)|^{1/2},
\end{equation}
the Eq. (\ref{pop}) is obtained, while the Eq. (\ref{pop1}) is reached for
\begin{equation}
\left|\det\dfrac{\partial\tilde{x}^\alpha(x)}{\partial x^\mu}\right|=|g(x)|.
\end{equation}
For this reason, one can induce a piecewise-Riemannian inner scalar product, norm and distance over $\mathcal{C}(\mathcal{M})$
\begin{subequations}
  \begin{align}
\langle X,Y\rangle&=\dfrac{\displaystyle\int\mathcal{D}X|Q(X)|^{1/2}\left|Q_{\mu\nu}X^\mu Y^\nu\right|}{\displaystyle\int\mathcal{D}X|Q(X)|^{1/2}},\\
||X||&=\left[\dfrac{\displaystyle\int\mathcal{D}X|Q(X)|^{1/2}\left|Q_{\mu\nu}X^\mu X^\nu\right|}{\displaystyle\int\mathcal{D}X|Q(X)|^{1/2}}\right]^{1/2}\label{norm},\\
d(X,Y)&=\left[\dfrac{\displaystyle\int\mathcal{D}X|Q(X)|^{1/2}\left|Q_{\mu\nu}(X^\mu-Y^\mu)(X^\mu-Y^\mu)\right|}{\displaystyle\int\mathcal{D}X|Q(X)|^{1/2}}\right]^{1/2}.
  \end{align}
\end{subequations}
Consequently, one can summarize
\newtheorem*{theorem*}{Theorem}
\begin{theorem*}
Every three non-singular charts on $\mathcal{M}$ generate a conformal space $\mathcal{C}(\mathcal{M})$ equipped with a diffeoinvariant Borel measure and a piecewise-Riemannian inner product.
\end{theorem*}

\section{Discussion}

Making use of the standard topological axioms, it can be shown directly that $\mathcal{C}(\mathcal{M})$ is a vector space, a function space, a real inner product space, a normed space, a real Hilbert space, and a metric space. By virtue of a piecewise-Riemannian nature, $\mathcal{C}(\mathcal{M})$ induced by the Borel measure is complete, as a result of the Cauchy--Schwarz inequality for inner product and the triangle inequality for distance.

Because $\mathcal{C}(\mathcal{M})$ is a normed metric space, it is a Hausdorff topological space with topology given by unions of open balls. Because all elements of a $\mathcal{C}(\mathcal{M})$ are the Cartesian products of the set of reals, $\mathcal{C}(\mathcal{M})$ is a finite topological space because the underlying set is finite. Because every finite topological space is a compact space, $\mathcal{C}(\mathcal{M})$ is a compact metric space, that is a metrizable space. Because a distance induces the Tychonoff topology, $\mathcal{C}(\mathcal{M})$ is completely metrizable space, that is a metrically topological complete space, and is homeomorphic to a complete metric space which has a countable dense subset.

Consequently, $\mathcal{C}(\mathcal{M})$ is a Polish space, a Luzin space, a Suslin space, and a Radon space. Among many other topological properties, it can be shown that $\mathcal{C}(\mathcal{M})$ is a separable Banach space, a perfectly normal Hausdorff space, a locally convex space, a bornological space, a reflexive space, a nuclear space, a barelled space, a perfect space, an Alexandroff space, a Baire space, a Borel space, a Dieudonn\'{e} complete space, a Fr\'{e}chet space, a Grothendieck space, a Mackey space, a Montel space, an Orlicz space, a Pt\'{a}k space, a Riesz space, a Schwartz space, a Sobolev space, an Urysohn space.

\end{document}